\documentclass[12pt,twoside]{article}
\usepackage{amsmath,amsbsy,amsfonts}
\usepackage{color}
\def\nd{\noindent}
\def\thend{\rule{3mm}{3mm}}
\def\Re{\mathbb{R}}

\def\Ka{\mathbb{K}}
\newtheorem{thm}{Theorem}[section]
\newtheorem{prop}{Proposition}[section]
\newtheorem{lem}{Lemma}[section]
\newtheorem{rmk}{Remark}[section]

\newtheorem{claim}{Claim}[section]
\newtheorem{cor}{Corollary}[section]
\newcommand{\fim}{\hfill\rule{2mm}{2mm}}
\newenvironment{dem}{\nd{\bf Proof. }}{\hskip.3cm\thend}

\begin{document}

\setlength{\baselineskip}{6.5mm} \setlength{\oddsidemargin}{8mm}
\setlength{\topmargin}{-3mm}
\title{\Large\sf Multiplicity of Positive Solutions for an \\ Obstacle Problem in $\Re$}
\author{\sf
Claudianor O. Alves\thanks{Partially supported by CNPq/Brazil
620150/2008-4 and 303080/2009-4} \, and Francisco Julio S.A. Corr\^ea \thanks{Partially supported by CNPq/Brazil  300561/2010-5 }\\
 Universidade Federal de Campina Grande\\
Centro de Ci\^encias e Tecnologia\\
Unidade Acad\^emica de Matem\'atica e Estat\'istica\\
 58.429-140 - Campina Grande - PB - Brazil\\
e-mail: coalves@dme.ufcg.edu.br, \,  fjsacorrea@gmail.com }

\pretolerance10000
\date{}
\numberwithin{equation}{section} \maketitle

\begin{center}
{\it Dedicated to Bernhard Ruf on the occasion of his 60th Birthday}
\end{center}

\begin{abstract}
In this paper we establish the existence of two positive solutions for the obstacle problem
$$
\displaystyle \int_{\Re}\left[u'(v-u)'+(1+\lambda V(x))u(v-u)\right] \geq \displaystyle \int_{\Re} f(u)(v-u), \forall v\in \Ka
$$
where $f$ is a continuous function verifying some technical conditions and $\Ka$ is the convex set given by
$$
\Ka =\left\{v\in H^{1}(\Re ); v \geq \varphi \right\},
$$
with $\varphi \in H^{1}(\Re)$ having nontrivial positive part with compact support in $\Re$.

\vspace{0.2cm} \noindent \emph{2000 Mathematics Subject
Classification} : 34B18, 35A15, 46E39.


\end{abstract}

\section{Introduction}

In this paper we will be concerned with the question of existence of positive solutions of a kind of obstacle problem. This class of problems has been largely studied due both its mathematical interest and its physical applications.  For example, it appears in mechanics, engineering, mathematical programming and optimization, among other things.  See, for instance, the classical books Kinderlehrer \& Stampacchia \cite{kinderlehrer}, Rodrigues \cite{rodrigues} and Troianiello \cite{troianiello} and the references therein.

 The typical obstacle problem is as follows: Let $\Omega$ be a domain in $\Re^{N}$. Given  functions $g: \Re \rightarrow \Re$ and $\varphi: \Omega \rightarrow \Re$, finding $u \in H_{0}^{1}(\Omega )$ satisfying
$$
\int_{\Omega}\nabla u \cdot \nabla (v-u)\geq \int_{\Omega}g(u)(v-u) \eqno{(P)}
$$
for all function $v$ in the convex set
\begin{equation} \label{E1.2}
\Ka := \left\{v\in H_{0}^{1}(\Omega ); v(x)\geq \varphi (x) \;\; a.e. \;\; \Omega \right\}
\end{equation}
where $\varphi$ is called the obstacle function.

Related to this kind of problem, the reader may consult Jianfu (\cite{jianfu1}, \cite{jianfu2}), where the author uses variational methods, Le \cite{Le} in which is used subsolution-supersolution techniques, Chang \cite{chang} where it is considered an obstacle problem related to discontinuous nonlinearities and Rodrigues \cite{rodrigues2} who considers combination of the obstacle problem with nonlocal equations in a class of free boundary problems. For more recent references we may cite Matzeu \& Servadei \cite{MS1}, in which the authors adapt for inequalities the iterative technique contained in de Figueiredo, Girardi \& Matzeu \cite{deFGM} for elliptic equations, Matzeu \& Servadei \cite{MS2} where the stability of solutions obtained in  \cite{MS1} are analized. Other results may be found in Servadei \& Valdinoci \cite{SV}, Mancini \& Musina \cite{MM}, Servadei (\cite{S}, \cite{ST}), Magrone, Mugnai \& Servadei \cite{MMS}.

These works and the references therein show clearly the mathematical importance and the wide variety of practical situations in which   obstacle problems may be found and applied.

Here we are interested in the unidimensional counterpart of problem $(P)$. More precisely, we consider the problem
$$
\displaystyle \int_{\Re}\left[u'(v-u)'+(1+\lambda V(x))u(v-u)\right] \geq \displaystyle \int_{\Re} f(u)(v-u), \forall v\in \Ka , \eqno{(P_\lambda)}
$$
where $u$ is a nonnegative function belonging to the convex set $\Ka$ given by
$$
\Ka :=\left\{v\in H^{1}(\Re ); v \geq \varphi \right\},  \eqno{(\Ka )}
$$
where $\varphi \in H^{1}(\Re )$ is assumed to have nontrivial positive part, that is,  $\varphi_{+}=\max \left\{\varphi \;, 0\right\}  \not\equiv   0$ . Moreover, $\lambda >0$ is a parameter and $f:\Re \to \Re$ is a nondecreasing continuous function verifying the following assumptions:
$$
\frac{f(t)}{t} \to 0\,\,\, \mbox{as} \,\,\, |t| \to 0  \eqno{(f_1)}
$$
and the Ambrosetti \& Rabinowitz Condition, that is, there is $\theta >2$ such that
$$
0<\theta F(t) \leq f(t)t \,\,\, \forall t \in \Re \setminus \{0\} \eqno{(f_2)}
$$
where $F(t)=\int_{0}^{t}f(s)ds$. We assume that  $V:\Re \rightarrow \Re$ is a nonnegative continuous function such that
$$
{\cal O} :=int \; \left(\left(V^{-1}(\left\{0\right\})\right)\right)\neq \emptyset
$$
is a bounded open set of $\mathbb{R}$ containing the support of $\varphi_{+}$, that is,  $Supp \; (\varphi_{+})\subset {\cal O}$. Here, $Supp (\varphi_{+})$ denotes the
support of $\varphi_{+}$ and
$$
V^{-1}(\left\{0\right\})=\left\{x\in \Re ; V(x)=0\right\}.
$$

The present paper was motivated by recent works involving the following class of problems
$$
\left\{
\begin{array}{l}
-\Delta{u}+(1+\lambda V(x))u=f(u) \,\,\, \mbox{in} \;\;\;\mathbb{R}^{N} \\
\mbox{}\\
u(x)>0 \,\,\, \mbox{in} \,\,\, \mathbb{R}^{N}
\end{array}
\right.
$$
where $\lambda $ is a positive parameter, $V:\mathbb{R}^{N} \to \mathbb{R}$ is a nonnegative function and $f$ is a continuous function satisfying some technical conditions. The reader may find more details in the papers of Alves \cite{A1}, Bartsch \& Wang \cite{BW2}, Clapp \& Ding \cite{CD}, Ding \& Tanaka \cite{DT} and their references. Here, we adapt some approaches found in these references to study the obstacle problem $(P_\lambda)$.

Our main result is the following

\begin{thm} \label{T1} Suppose $(f_1)-(f_2)$ hold, then there are $r, \lambda^{*} >0$, such that if  $\|\varphi_{+} \|_{H^{1}(\Re)} < \frac{r}{2}$, problem $(P_{\lambda})$ has two positive solutions for all $\lambda > \lambda^{*}$ .

\end{thm}

One of the main difficulties to prove Theorem \ref{T1} is related to the fact that the energy functional associated with the problem $(P_{\lambda})$ does not satisfy in general the well known Palais-Smale condition, once that we are working in whole $\mathbb{R}$. To overcome this difficulty,  we adapt some ideas found in {del Pino \& Felmer \cite{DF}}, modifying the function $f$ outside the set ${\cal O}$, in such way that the energy functional of the modified obstacle problem satisfies the Palais-Smale condition.  Using variational methods, we prove the existence of two solutions for the modified obstacle  problem. After that, it is proved that under the hypotheses of Theorem \ref{T1}, the solutions found are solutions of the original obstacle problem.

The structure of this paper is as follows: In Section 2 we introduce the modified obstacle problem, in Section 3 we establish the existence of a first solution for the modified obstacle problem by minimization, in Section 4 we show the existence of a second solution for the modified obstacle problem by the  Mountain Pass Theorem and in Section 5 we prove Theorem \ref{T1}.

\section{The  Modified Obstacle Problem}

From this time onwards, since we intend to find positive solution, we will assume, without loss of generality, that
$$
f(t)=0 \,\,\, \forall t \leq 0.
$$

To prove the existence of positive solutions for $(P_{\lambda})$, we will work with a modified obstacle problem, following some ideas found in del Pino \& Felmer \cite{DF}. To this end, we consider the function $h:\Re \rightarrow \Re$ as follows:
$$
h(t)=\left\{
\begin{array}{ccc}
f(t) & \mbox{if} & t\leq a,\\
\frac{1}{k}t & \mbox{if} & t\geq a,
\end{array}
\right.
$$
where $k > \max\{\frac{\theta}{\theta -2},2\}$ and $a >0$  satisfy  $\frac{f(a)}{a}=\frac{1}{k}$. We now set
$$
g(x,t)=\chi(x)f(t)+(1-\chi(x))h(t),
$$
where $\Omega \subset \mathbb{R}$ is a bounded open set containing ${\cal O}$ and $\chi$ is the characteristic function of the set $\Omega$, that is,
$$
\chi(x)=
\left\{
\begin{array}{l}
1, \;\;\; x \in \Omega \\
0, \;\;\; x \in \Omega^{c}.
\end{array}
\right.
$$

Using the function $g$, we will show the existence of two positive solutions for the obstacle problem
$$
\displaystyle\int_{\Re}\left[u'(v-u)'+(1+\lambda V(x))u(v-u)\right] \geq \displaystyle\int_{\Re}g(x,u)(v-u),  \,\,\, \forall v \in \Ka . \leqno{(P_{A})}
$$

\begin{rmk} \label{R1}
If $u $ is a solution of $(P_{A})$ verifying
$$
u(x)\leq a, \; \forall x \in \Omega^{c},
$$
then $u$ is a solution of the original obstacle problem. Indeed,
if $x\in \Omega$, we have $\chi(x)=1$ and so
$$
g(x,u(x))=f(u(x)).
$$
If $x\notin \Omega \; (x\in \Omega^{c})$, then $\chi(x)=0$ and so
$$
g(x,u(x))=h(u(x))=f(u(x)),
$$
because $h(u(x))=f(u(x))$ since $0\leq u(x)\leq a$ in $\Omega^{c}$.
\end{rmk}

Let $E_\lambda \subset H^{1}(\Re )$ be the subspace
$$
E_{\lambda}=\left\{u\in H^{1}(\Re ); \displaystyle \int_{\Re}V(x)u^{2}<\infty \right\}
$$
endowed with the norm
$$
\|u\|_{\lambda}=\Big( \int_{\mathbb{R}}[|u'|^{2}+ (1+\lambda V(x))|u|^{2}] \Big)^{\frac{1}{2}}.
$$
Hereafter, we denote by $\| \,\,\, \|$ the usual norm in $H^{1}(\Re)$.

Since we approach our problem by means of variational method, we consider the energy functional associated with the obstacle problem $(P_{A})$, $I_{\lambda}:E_{\lambda} \rightarrow \Re$,  given by
$$
I_{\lambda}(u)=\frac{1}{2}\|u\|_{\lambda}^{2}-\displaystyle \int_{\Re}G(x,u) + \Psi (u),
$$
where
$$
G(x,t)=\displaystyle \int_{0}^{t}g(x,s)ds
$$
and $\Psi : E\rightarrow (-\infty , \infty ]$ is the indicatrix function of the set $\Ka$, i.e.,
\begin{equation} \label{indicadora}
\Psi (u)=0, \; \forall u \in \Ka \,\,\, \mbox{and} \,\,\, \Psi (u)=+\infty , \; \forall u \in \Ka^{c}.
\end{equation}

\begin{prop} \label{PS}
The functional $I_\lambda$ satisfies the $(PS)$ condition.
\end{prop}

\dem Let $d \in \mathbb{R}$ and $(u_{n}) \subset E_\lambda$ be a $(PS)_d$ sequence for $I_\lambda$. Then, there is $(z_{n})\subset E_\lambda'$ with $z_n \to 0$ such that
$$
I_\lambda(u_{n})\rightarrow d \,\,\, \mbox{and} \,\,\, I_\lambda'(u_n)(v-u_n) \geq <z_{n},v-u_{n}> \,\,\, \forall n \in \mathbb{N} \,\,\,\, \mbox{and} \,\,\, v \in \Ka ,
$$
that is,
$$
\displaystyle \int_{\Re}u_{n}'(v-u_{n})'+(1+\lambda V(x))u_{n}(v-u_{n})-\displaystyle \int_{\Re}g(x,u_{n})(v-u_{n})\geq <z_{n},v-u_{n}>,
$$
for all $v \in \Ka$.

\begin{claim} \label{C1}
$(u_{n})$ is a bounded sequence in $E_\lambda$.
\end{claim}

\noindent We deal separately with the sequences $({u_{n}}_{+})$ and $({u_{n}}_{-})$, where ${u_{n}}_{-}=\max\{-u_n,0\}$ . Since $u_{n}={u_{n}}_{+}-{u_{n}}_{-}$, it is enough to show that $({u_{n}}_{+})$ and $({u_{n}}_{-})$ are bounded in $E_\lambda$. To show the boundedness of  $({u_{n}}_{-})$, we consider the test function $v=u_{n}+{u_{n}}_{-} \in \Ka$. So,
$$
\displaystyle \int_{\Re}(u_{n}'({u_{n}}_{-})'+(1+\lambda V(x))u_{n}{u_{n}}_{-})-\int_{\Re}g(x,u_{n}){u_{n}}_{-}\geq <z_{n},{u_{n}}_{-}>.
$$
Because $\displaystyle \int_{\Re}(1+\lambda V(x)){u_{n}}_{+}{u_{n}}_{-}=\displaystyle \int_{\Re}g(x,u_{n})u_{n}^{-}=0$, we obtain
$$
-\|{u_{n}}_{-}\|_{\lambda}^{2}\geq <z_{n},{u_{n}}_{-}>,
$$
which leads to
$$
\|{u_{n}}_{-}\|_{\lambda}^{2}\leq \|z_{n}\|\|{u_{n}}_{-}\|_{\lambda}.
$$
As $z_{n}\rightarrow 0$ in $E_\lambda'$, we conclude that ${u_{n}}_{-}\rightarrow 0$ in $E_{\lambda}$, and thus, $({u_{n}}_{-})$ is bounded in $E_\lambda$.

With respect to $({u_{n}}_{+})$, fixing the test function $v=u_{n}+{u_{n}}_{+} \in \Ka$, we derive that
\begin{equation} \label{E10}
\|{u_{n}}_{+}\|_{\lambda}^{2}-\displaystyle \int_{\Re}g(x,u_{n})u_{n}^{+}\geq <z_{n}, u_{n}^{+}>,
\end{equation}
leading to
\begin{equation}\label{E11}
- \int_{\Omega}f(u_n)u_n \geq - \|{u_{n}}_{+}\|_{\lambda}^{2} + <z_{n}, u_{n}^{+}>.
\end{equation}
On the other hand, we know that
$$
d=\frac{1}{2}\|u_{n}\|_{\lambda}^{2}-\int_{\Omega}F(u_n)-\int_{\Omega^{c}}G(x,u_{n}) +o_{n}(1).
$$
Using the definition of $g$, it is easy to prove that
\begin{equation} \label{G}
2G(x,t) \leq g(x,t)t \leq \frac{1}{k}(1+\lambda V(x))|t|^{2} \,\,\, \forall x \in \Omega^{c} \,\,\, \mbox{and} \,\,\, t \in \Re.
\end{equation}
Thereby, from $(f_{2})$ and (\ref{G})
\begin{equation} \label{E12}
d \geq \frac{1}{2}\|{u_{n}}_{+}\|_{\lambda}^{2}-\frac{1}{\theta}\int_{\Omega}f(u_n)u_n -\frac{1}{2k}\int_{\Omega^{c}}(1+\lambda V(x))|u_{n}|^{2} +o_{n}(1).
\end{equation}
Combining (\ref{E11}) and (\ref{E12}),
$$
d \geq \left[ \Big( \frac{1}{2} - \frac{1}{\theta}\Big) - \frac{1}{2k} \right]  \|{u_{n}}_{+}\|_{\lambda}^{2} - \|z_{n}\| \|{u_{n}}_{+}\| + o_{n}(1).
$$
Since $ k > \frac{\theta}{\theta - 2}$ and $z_n \to 0$ in $E_{\lambda}'$, the last inequality implies that $({u_n}_{+})$ is bounded in $E_{\lambda}$. Therefore, $(u_n)$ is  bounded in $E_{\lambda}$.

Now, we will show that $(u_n)$ has a subsequence that converges strongly in $E_{\lambda}$.  {Since  $({u_{n}}_{-})$ converges to 0 in $E_\lambda$, without loss of generality, we will assume that $u_n \geq 0$ for all $n \in \mathbb{N}$}. We begin by fixing $R>0$ so large in order that
$\Omega \subset \left(-\frac{R}{2}, \frac{R}{2}\right)$ and a function $\eta \in C^{1}(\Re, \Re)$ satisfying

\vspace{0.5cm}

$\bullet$ $0\leq \eta (t) \leq 1, \forall t\in \Re ;$

\vspace{0.5cm}

$\bullet$  $\eta (t)=0, |t|\leq \frac{R}{2};$

\vspace{0.5cm}

$\bullet$ $\eta (t)=1, |t|\geq R ;$

\vspace{0.5cm}

$\bullet$ $ |\eta' (t)|\leq \frac{C}{R}, \forall t\in \Re $.

\begin{claim} \label{D1} Given $\delta >0$, there is $R>0$ such that
$$
\displaystyle \int_{|x|\geq {R}}(|u_{n}'|^{2}+|u_{n}|^{2})<\delta .
$$
\end{claim}
Assuming that this claim is true, we continue with our proof. Considering the test function $v= u_n-\eta (u_n - \varphi_{+}) =u_n - \eta u_n \in \Ka$, it follows that
$$
\displaystyle \int_{\Re}\left[u_{n}'(\eta u_{n})'+(1+\lambda V(x))u_{n}(\eta u_{n})\right]\leq \int_{\Re}g(x,u_{n})(\eta u_{n})+o_{n}(1)
$$
or, equivalently,
$$
\int_{\Re}\eta |u_{n}'|^ {2}+\int_{\Re}u_{n}'\eta'u_{n}+\int_{\Re}(1+\lambda V(x))\eta |u_{n}|^ {2}\leq \int_{|x|\geq \frac{R}{2}}g(x,u_{n})\eta u_{n} +o_{n}(1)
$$
implying that
$$
\begin{array}{l}
\displaystyle \int_{|x|\geq R}|u'_{n}|^{2}+\int_{|t|\leq R}u_{n}'\eta' u_{n}+ \int_{|x|\geq \frac{R}{2}}(1+\lambda V(x))\eta |u_{n}|^{2}\leq \\
\mbox{}\\
\displaystyle \int_{|x|\geq \frac{R}{2}}\frac{1}{k}(1+\lambda V(x)) \eta |u_{n}|^{2} +o_{n}(1).
\end{array}
$$
Because $k >2$, it follows that
$$
\begin{array}{l}
\displaystyle \int_{|x|\geq R}|u_{n}'|^{2}+\int_{|t|\leq R}u_{n}'\eta' u_{n}+\int_{|x|\geq \frac{R}{2}}(1+\lambda V(x))\eta |u_{n}|^{2}\leq \\
\mbox{}\\
\displaystyle \int_{|x|\geq \frac{R}{2}}\left(\frac{1+\lambda V(x)}{2}\right)|u_{n}|^{2}+o_{n}(1)
\end{array}
$$
and so,
$$
\int_{|x|\geq R}|u_{n}'|^{2}+\frac{1}{2}\int_{|x|\geq \frac{R}{2}}(1+\lambda V(x))\eta |u_{n}|^{2}\leq \int_{|x|\leq R}|u_{n}'||\eta'||u_{n}|\leq \frac{C}{R}+o_{n}(1).
$$
Thereby,
$$
\int_{|x|\geq R}|u_{n}'|^{2}+\int_{|x|\geq R}(1+\lambda V(x))|u_{n}|^{2} \leq \frac{{C}}{R}+o_{n}(1),
$$
showing that
$$
\limsup_{n \to +\infty} \int_{|x|\geq R}(|u_{n}'|^ {2}+|u_{n}|^ {2})\leq \frac{{C}}{R}.
$$
Now, we choose $R>0$ so large in order
$$
\limsup_{n \to +\infty} \int_{|x|\geq R}(|u_{n}'|^ {2}+|u_{n}|^ {2})<\delta ,
$$
{proving the Claim \ref{D1}.}

Recalling that for each $R>0$, the Sobolev embedding
$$
H^{1}(\mathbb{R}) \hookrightarrow C([-R,R])
$$
is compact, we have that
$$
u_{n}\rightarrow u \,\,\, \mbox{in} \,\,\, C([-R,R]).
$$
{This limit, combined with the Claim \ref{D1}}, asserts that
\begin{equation} \label{E130}
\int_{\Re}g(x,u_{n})u_{n}\rightarrow \int_{\Re}g(x,u)u
\end{equation}
and
\begin{equation} \label{E13}
\int_{\Re}g(x,u_{n})v\rightarrow \int_{\Re}g(x,u)v \,\,\, \forall v \in \Ka,
\end{equation}
where $u \in \Ka$ is the weak limit of $(u_n)$ in $E_{\lambda}$.

Since $(u_n)$ is a bounded Palais-Smale sequence for $I_{\lambda}$, we have
\begin{equation} \label{E14}
\int_{\Re}u_{n}'(v-u_{n})'+(1+\lambda V(x))u_{n}(v-u_{n})\geq \int_{\Re}g(x,u_{n})(v-u_{n})+o_{n}(1) \,\,\, \forall  v \in \Ka
\end{equation}
or equivalently
$$
\begin{array}{l}
\displaystyle \int_{\Re}[u_{n}'v'+(1+\lambda V(x))u_{n}v] \geq \int_{\Re}[|u_n'|^{2} + (1+\lambda V(x))|u_{n}|^{2}] \, + \\
\mbox{}\\
\displaystyle + \int_{\Re}g(x,u_{n})(v-u_{n})+o_{n}(1)
\end{array}
$$
for all $v \in \Ka.$ Taking the inferior limits on both sides of the above inequality and using (\ref{E130}) and (\ref{E13}),  we get
$$
\begin{array}{l}
\displaystyle \int_{\Re}[u'v'+(1+\lambda V(x))uv] \geq \int_{\Re}[|u'|^{2} + (1+\lambda V(x))|u|^{2}] \, + \\
\mbox{}\\
\displaystyle + \int_{\Re}g(x,u)(v-u)+o_{n}(1)
\end{array}
$$
that is,
$$
\int_{\Re}[u'(v-u)'+(1+\lambda V(x))u(v-u)]\geq \int_{\Re}g(x,u)(v-u), \,\,\, \forall v \in \Ka
$$
from where it follows that $u$ is a critical point of $I_{\lambda}$.

Using  $u$ as a test function in (\ref{E14}) and the limit (\ref{E13}), it follows that
$$
\limsup_{n \to +\infty} \|u_n\|_{\lambda}^{2} \leq \|u\|_{\lambda}^{2}.
$$
Since $E_{\lambda}$ is a Hilbert space, the last inequality leads to $u_n \to u$ in $E_{\lambda}$, finishing the proof of proposition.  {\fim}

\section{First solution for $(P_{A})$}

The first positive solution of $(P_{A})$ will be obtained via Ekeland's Variational Principle \cite{Ekeland}.  In this section, we denote by $B_{r}$ and $\Ka_{r}$ the following sets
$$
B_{r}=\{u \in E_{\lambda}; \, \|u\|_{\lambda} < r \}
$$
and
$$
\Ka_{r} = \Ka \cap \overline{B}_{r} .
$$

\begin{thm}  \label{primeirasolucao} There is $r>0$, such that if $\|\varphi_{+}\|_{H^{1}(\Re)}<\frac{1}{2}\sqrt{r}$, the variational problem
\begin{equation} \label{Minimizacao}
m=\inf\{I_{\lambda}(u): \, u \in \Ka_{r}\}
\end{equation}
has a solution for all $\lambda >0$. Moreover, this solution is a positive solution of $(P_A)$.
\end{thm}

\dem  First of all, we observe that
$$
\int_{\Re}G(x,u(x))=\displaystyle \int_{\Omega}F(u)+\displaystyle \int_{\Omega^{c}}G(x,u(x)).
$$
From $(f_1)$, if $\|u\|_{\lambda}=r$ and $r$ is small enough, we have that
$$
\int_{\Omega}F(u) \leq \frac{1}{4}\int_{\Omega}|u|^{2} \leq \frac{1}{4}\|u\|^{2}_{\lambda}.
$$
Hence
$$
\int_{\Re}G(x,u(x)) \leq \frac{1}{4}\|u\|^{2}_{\lambda} + \int_{\Omega^{c}}G(x,u(x)),
$$
and so, by (\ref{G}),
$$
\int_{\Re}G(x,u(x)) \leq \frac{1}{4}\|u\|^{2}_{\lambda}+ \frac{1}{2k}\int_{\Omega^{c}}(1+\lambda V(x))|u|^{2}.
$$
Thereby,
\begin{equation}
I_\lambda(u)\geq \frac{1}{4}\|u\|_{\lambda}^{2}-\frac{1}{2k}\int_{\Omega^{c}}(1+\lambda V(x))|u|^{2}+\Psi (u)
\end{equation}
from where it follows that
\begin{equation}
I_\lambda(u)\geq \Big( \frac{1}{4} - \frac{1}{2k} \Big) \|u\|_{\lambda}^{2}+\Psi (u), \forall u\in E_{\lambda}.
\end{equation}
Since $k > 2$,
\begin{equation} \label{AA1}
I(u)\geq \frac{1}{8}\|u\|_{\lambda}^{2}, \forall u\in \Ka_{r}.
\end{equation}

From the above study, we have that $m$ is well defined and {$m \in [0,+\infty)$}. Therefore, there is $(u_n) \subset \Ka_{r}$ such that
$$
I_{\lambda}(u_n) \to m.
$$
Once that $(u_n)$ is bounded, {because $(u_n) \subset \overline{B}_r(0)$}, we can assume, without loss of generality, that
$$
u_n \rightharpoonup u \,\,\, \mbox{in} \,\,\, E_{\lambda}
$$
and
$$
u_n(x) \to u(x) \,\,\, \mbox{a.e. in} \,\,\, \Re.
$$
By Ekeland's Variational Principle, we also assume that
$$
m \leq I_{\lambda}(u_n) \leq m + \frac{1}{n^{2}} \,\,\, \forall n \in \mathbb{N}
$$
and
$$
I_\lambda(u)\geq I_{\lambda}(u_n) - \frac{1}{n}\|u-u_n\|_{\lambda} \,\,\, \forall u \in \Ka_{r}.
$$
Observing that $\varphi_{+} \in \Ka_{r}$, by (\ref{AA1}),
$$
\frac{1}{8}\|u_n\|_{\lambda}^{2} \leq I_{\lambda}(u_n) \leq m +\frac{1}{n^{2}}\leq I_{\lambda}(\varphi_{+})+\frac{1}{n^{2}} \leq \frac{1}{2}\|\varphi_{+}\|^{2} +\frac{1}{n^{2}}
$$
leading to
$$
\limsup_{n \to +\infty} \|u_n\|_{\lambda}^{2} \leq 4 \, \|\varphi_{+}\|^{2}< r .
$$
Thus, there is $n_0 \in \mathbb{N}$ such that 
$$
 \|u_n\|_{\lambda}^{2}  < r \,\,\  \forall n \geq n_0.
$$
Now, repeating the same arguments found in \cite{jianfu2}, we have that $(u_n)$ is a $(PS)_m$ sequence for $I_{\lambda}$, that is,
\begin{equation} \label{F1}
I_\lambda(u_n) \to m \,\,\, \mbox{and} \,\,\,  I'_{\lambda}(u_n)(v-u_n) \geq <z_n, v-u_n> \,\,\, \forall v \in \Ka
\end{equation}
with $z_n \to 0$ in $E'_{\lambda}$. Using Proposition \ref{PS}, there are a subsequence of $(u_n)$, still denoted by $(u_n)$, and $u$ in $E_{\lambda}$ such that
\begin{equation} \label{F2}
u_n \to u \,\,\, \mbox{in} \,\,\, E_{\lambda}.
\end{equation}
From this, $u \in \Ka_{r}$ and $I_{\lambda}(u)=m$, showing that $u$ is a solution for (\ref{Minimizacao}).  Now, combining (\ref{F1}) and (\ref{F2}), it follows that
\begin{equation} \label{Sol1}
\int_{\Re}[u'(v-u)' + (1+\lambda V(x))u(v-u)]  \geq \int_{\Re}g(x,u)(v-u) \,\,\, \forall v \in \Ka.
\end{equation}
Using the test function $v=u +u_{-} \in \Ka$, a direct computation implies that $u_{-}=0$, consequently $u$ is nonnegative.  The positivity of $u$ is obtained by applying the maximum principle. \fim

\section{Second solution for $(P_{A})$}

In this section, we will apply the Mountain Pass Theorem due to Szulkin \cite{Szulkin} to get a second positive solution for problem $(P_A)$. Here, we denote by $u_\lambda$ the solution obtained in Theorem \ref{primeirasolucao}.

\begin{lem} \label{geometria}
The energy functional $I_{\lambda}$ verifies the geometry of the Mountain Pass Theorem.
\end{lem}

\dem {Note that, by Theorem \ref{primeirasolucao},
$$
I_{\lambda}(u) \geq I_{\lambda}(u_{\lambda}) \,\,\, \forall u \in \Ka_{r}.
$$
Since $\Psi(u)=+\infty$ for all $u \in \Ka_{r}^{c}$, it follows that
\begin{equation} \label{AA2}
I_{\lambda}(u) \geq I_{\lambda}(u_{\lambda}) \,\,\, \forall u \in \overline{B}_r.
\end{equation}
Moreover, if $\rho=\frac{1}{8}r^{2}$,  (\ref{AA1}) gives
$$
I_\lambda(u)\geq \rho >0, \;\;  \;\;\; \mbox{for all } \;\;\; u \in \partial \overline{B}_r.
$$
On the other hand, since $\|\varphi_{+}\|^{2}< \frac{1}{4}r^{2}$, we have that $\varphi_{+} \in \Ka_{r}$, and so,
\begin{equation} \label{Z0}
I_{\lambda}(u_{\lambda}) \leq I_{\lambda}(\varphi_{+})\leq \frac{1}{2}\|\varphi_{+}\|^{2}< \rho,
\end{equation}
from where it follows that
\begin{equation} \label{Z2}
\inf_{u \in \partial B_{r}}I_{\lambda}(u) > I_{\lambda}(u_{\lambda}).
\end{equation}
We now observe that, for $t\geq 1$, $t\varphi_{+}\in \Ka$. Then, $\Psi (t\varphi_{+})=0$ and
$$
I_\lambda(t\varphi_{+})=\frac{t^{2}}{2}\int_{\Re}(|\varphi_{+}'|^{2}+|\varphi_{+}|^{2})- \int_{\Re} F(t\varphi_{+}).
$$
By $(f_2)$, there are $A,B>0$ such that
$$
F(s) \geq As^{\theta} - B \,\,\, \forall s \geq 0.
$$
Consequently,
$$
I_\lambda(t\varphi_{+}) \leq \frac{t^{2}}{2}\int_{\Re}(|\varphi_{+}'|^{2}+|\varphi_{+}|^{2})- t^{\theta} A \int_{D} (\varphi_{+})^{\theta} + B|D|,
$$
where $D$ is a mensurable set with finite measure verifying $D \cap Supp (\varphi_{+}) \not= \emptyset$. From this,
$$
I_\lambda(t\varphi_{+})\rightarrow -\infty \;\; \mbox{as}\;\; t\rightarrow +\infty ,
$$
and thus, setting $e=t \varphi_{+}$ for $t$ large enough, we derive that
\begin{equation} \label{Z1}
\|e\| >r \,\,\, \mbox{and} \,\,\, I_{\lambda}(e) < I_{\lambda}(u_{\lambda}).
\end{equation}
From (\ref{AA2})-(\ref{Z1}), we deduce that $I_{\lambda}$ satisfies the mountain pass geometry, see \cite[Theorem 3.2]{Szulkin}. \fim }

\begin{thm} \label{T2} Under the assumptions of Theorem \ref{primeirasolucao}, Problem $(P_A)$ has a  positive solution at the mountain pass level for all $\lambda >0$, that is, there is $w_{\lambda} \in \Ka$ verifying
$$
I_{\lambda}(w_{\lambda})=c_{\lambda} \,\,\, \mbox{and} \,\,\, \,\,\, I'_{\lambda}(w_{\lambda})(v-w_\lambda) \geq 0 \,\,\, \forall v \in \Ka,
$$
where $c_{\lambda}$ is the mountain pass level of $I_{\lambda}$.
\end{thm}

\dem  Combining Lemma  \ref{geometria} and Propostion \ref{PS} with the Mountain Pass Theorem, we have that the mountain pass level $c_{\lambda}$ associated with $I_{\lambda}$ is a critical value, hence there is $w_{\lambda} \in \Ka$ such that
$$
I_{\lambda}(w_{\lambda})=c_{\lambda} \,\, \mbox{and} \,\,\, I'_{\lambda}(w_{\lambda})(v-w_\lambda) \geq 0 \,\,\, \forall v \in \Ka.
$$
Using the test function $v=w_\lambda +{w_{\lambda}}_{-} \in \Ka$, a direct computation implies that ${w_{\lambda}}_{-}=0$, consequently $w_\lambda$ is nonnegative.  The positivity of $w_\lambda$ is obtained by applying maximum principle. \fim

\begin{cor} \label{2solucoes} Under the assumptions of Theorem \ref{primeirasolucao}, problem $(P_A)$ has two positive solutions for all $\lambda >0$.
\end{cor}

\dem From the previous study, we have two solutions denoted by $u_\lambda$ and $w_\lambda$, where $u_{\lambda}$ was obtained by minimization and $w_{\lambda}$ by Mountain Pass Theorem. Moreover, {by (\ref{Z0}),}
$$
m=I_\lambda (u_\lambda) < \rho
$$
and
$$
I_\lambda (w_\lambda)=c_\lambda \geq \rho.
$$
Thus,
$$
I_\lambda(u_\lambda) < I_\lambda (w_\lambda),
$$
from where it follows that $u_\lambda$ and $w_\lambda$ are different. Hence, problem $(P_A)$ has two positive solutions. \fim

\section{Proof of Theorem \ref{T1}}

In what follows,  our main goal is to show that there is $\lambda^{*} >0$ such that if $\lambda \geq \lambda^{*}$, the solutions $u_{\lambda}$ and $w_{\lambda}$ obtained in Corollary \ref{2solucoes} satisfy
\begin{equation} \label{Condicao}
w_\lambda(x), u_{\lambda}(x)\leq a, \; \forall x \in \Omega^{c}.
\end{equation}
From this, by using  Remark \ref{R1}, we will conclude that $w_\lambda$ and $u_{\lambda}$ are positive solutions of $(P_{\lambda})$ if $\lambda \geq \lambda^{*}$.

Hereafter, $\lambda_{n} \to +\infty$, $u_{n}=u_{\lambda_{n}}$ and $w_{n}=w_{\lambda_{n}}$. From Theorem \ref{primeirasolucao}, we know that $u_n \in \Ka_r$ for all $n \in \mathbb{N}$, thus $(u_n)$ is bounded in $H^{1}(\Re)$. Next, we will show that $(w_n)$ is also bounded in $H^{1}(\Re)$.

\begin{lem}
\label{limita}
The sequence $(w_n)$ is bounded in $H^{1}(\Re)$. More precisely, there is $M>0$ such that
$$
\|w_n\|_{\lambda_{n}} \leq M \,\,\, \forall n \in \mathbb{N}.$$

\end{lem}

\dem Since $w_n$ is a solution of $(P_{\lambda_{n}})$, it follows that
\begin{equation} \label{X1}
\int_{\Re}[w_n'(v-w_n)'+(1+\lambda_n V(x))w_n(v-w_n)]\geq \int_{\Re}g(x,w_n)(v-w_n), \,\,\, \forall v \in \Ka.
\end{equation}
Repeating the same arguments used in the proof of Proposition \ref{PS}, we derive that
\begin{equation} \label{E15}
I_{\lambda_{n}}(w_n) \geq \left[ \Big( \frac{1}{2} - \frac{1}{\theta}\Big) - \frac{1}{2k} \right]  \|w_{n}\|_{\lambda_{n}}^{2} \,\,\, \forall n \in \mathbb{N}.
\end{equation}
Now, considering the path $\gamma (t)= t t^{*} \varphi_{+}$ for $t \in [0,1]$ and $t^{*}$ large enough  and setting
$$
\Sigma = \max_{t \in [0,1]} J(\gamma (t)) >0, 
$$
where
$$
J(u)= \frac{1}{2} \int_{\Omega} [|u'|^{2}+|u|^{2}] - \int_{\Omega}F(u),
$$
{it follows that
$$
I_{\lambda_{n}}(w_n) \leq \max_{t \in [0,1]} I_{\lambda_n}(\gamma (t))=\max_{t \in [0,1]} J(\gamma (t))= \Sigma \,\,\, \forall n \in \mathbb{N},
$$
because $I_{\lambda_n}(\gamma(t))=J(\gamma (t))$ for all $n \in \mathbb{N}$ and $t \in [0,1]$ .}

This combined with (\ref{E15}) implies that $(\|w_n\|_{\lambda_{n}})$ is bounded in $\Re$. \fim

\begin{lem} \label{convergencia} There are subsequences of $(u_n)$ and $(w_n)$, still denoted by itself, which are strongly convergent in $H^{1}(\Re)$.
\end{lem}

\dem In what follows, we will prove the lemma only for $(u_n)$, because the same arguments can be applied to $(w_n)$. Following the same arguments used in the proof of Proposition  \ref{PS}, for each $\delta >0$, there is $R>0$ such that
$$
\limsup_{n \to +\infty} \int_{|x|\geq R}[|u_n'|^{2}+|u_n|^{2}] < \delta.
$$
The above limit yields
\begin{equation} \label{E16}
\int_{\Re}g(x,u_{n})u_{n}\rightarrow \int_{\Re}g(x,u)u
\end{equation}
and
\begin{equation} \label{E17}
\int_{\Re}g(x,u_{n})v\rightarrow \int_{\Re}g(x,u)v \,\,\, \forall v \in \Ka,
\end{equation}
where $u \in \Ka$ is the weak limit of $(u_n)$ in $H^{1}(\Re)$.

\begin{claim} \label{C2}
The weak limit $u$ is null in $\overline{{\cal O}^{c}}$, that is,
$$
u(t)=0 \,\,\, \forall t \in \overline{{\cal O}^{c}}.
$$
Hence, $u \in H^{1}_{0}({\cal O})$.
\end{claim}

\noindent In fact, for each $m \in \mathbb{N}$, we define
$$
\Delta_{m}=\left\{t \in \Re; \, V(t) > \frac{1}{m} \right\}.
$$
It is immediate to see that
$$
P=\{t \in \Re; \, V(t) >0\}= \displaystyle \bigcup_{m=1}^{\infty} \Delta_{m}.
$$
Note that
$$
\int_{\Delta_{m}}|u_n|^{2}\leq \frac{m}{\lambda_{n}}\|u_n\|_{\lambda_{n}}^{2}\leq \frac{m}{\lambda_{n}}r^{2} \,\, \forall n,m \in \mathbb{N}
$$
where {$r$ is the constant given in Theorem \ref{primeirasolucao}}. The last inequality, together with Fatou's Lemma, lead to
$$
\int_{\Delta_{m}}|u|^{2}=0 \,\,\, \forall m \in \mathbb{N}.
$$
Thereby, $u=0$ a.e in $\Delta_{m}$ for all $m \in \mathbb{N}$, implying that $u=0$ a.e. in $P$. Now, the claim follows using the continuity of $u$.

Using $v=u$ as a test function in (\ref{Sol1}),
\begin{equation} \label{Z3}
\int_{\Re}|u_{n}'|^{2}+ \int_{\Re}(1+\lambda_{n} V)|u_{n}|^{2} \leq  \int_{\Re}(1+\lambda_n V)u_{n}u + \int_{\Re}u_{n}'u'  - \int_{\Re}g(x,u_{n})(u-u_{n}).
\end{equation}
Once that $V(t) \geq 0$ and $u=0$ in $\overline{\Omega^{c}}$,
$$
\int_{\Re}|u_{n}'|^{2}+ \int_{\Re}|u_{n}|^{2} \leq \int_{\Re}u_{n}'u' +\int_{\Re}u_{n}u  - \int_{\Re}g(x,u_{n})(u-u_{n}).
$$
Taking the limit of $n \to +\infty$ and using (\ref{E16})-(\ref{Z3}),
$$
\limsup_{n \to +\infty} \int_{\Re}[|u_{n}'|^{2}+|u_{n}|^{2}] \leq  \int_{\Re}[|u'|^{2}+ |u|^{2}].
$$
Since $H^{1}(\Re)$ is a Hilbert space and $u_n \rightharpoonup u$ in $H^{1}(\Re)$, the above limit implies that $u_n \to u$ in $H^{1}(\Re)$. \fim

As a consequence of the lemmas proved in this section, we have the following results

\begin{cor} \label{Corolario1} {The sequences $(u_n)$ and $(w_n)$ satisfy
\begin{equation} \label{W1}
\lambda_{n}\int_{\Re}V(x)|u_n|^{2} \to 0 \,\,\, \mbox{as} \,\,\, n \to +\infty
\end{equation}
and
\begin{equation} \label{Ww1}
\lambda_{n}\int_{\Re}V(x)|w_n|^{2} \to 0 \,\,\, \mbox{as} \,\,\, n \to +\infty,
\end{equation}
for some subsequence. Moreover, the weak limits $u$ and $w$ of $(u_n)$ and $(w_n)$ respectively, belong to $H^{1}_{0}({\cal O})$ and they are positive solutions of the obstacle problem
$$
\int_{{\cal O}}[\psi'(v-\psi)' + \psi(v-\psi)] \geq \int_{{\cal O}}f(\psi)(v-\psi) \,\,\, \forall v \in \widehat{\Ka} \eqno{(P_{{\cal O}})}
$$
where
$$
\widehat{\Ka} := \left\{v\in H_{0}^{1}({\cal O} ); v(x)\geq \varphi (x) \;\; a.e. \;\; {\cal O} \right\}.
$$}
\end{cor}

\dem From now on, we will prove the lemma only for the sequence $(u_n)$, because the same arguments can be applied to $(w_n)$. {Repeating the same type of arguments explored in the proof of Claim \ref{C2}, we get again an equality like (\ref{Z3}), that is,
$$
\int_{\Re}|u_{n}'|^{2}+ \int_{\Re}(1+\lambda_{n} V)|u_{n}|^{2} \leq  \int_{\Re}(1+\lambda_n V)u_{n}u + \int_{\Re}u_{n}'u'  - \int_{\Re}g(x,u_{n})(u-u_{n}).
$$
Using the fact that $V(t)u(t)=0$ for all $t \in \mathbb{R}$, it follows that
\begin{equation} \label{Z4}
\int_{\Re}|u_{n}'|^{2}+ \int_{\Re}(1+\lambda_{n} V)|u_{n}|^{2} \leq \int_{\Re}u_{n}'u' + \int_{\Re}u_{n}u  - \int_{\Re}g(x,u_{n})(u-u_{n}).
\end{equation}
From Theorem \ref{convergencia}, $u_n \to u$ in $H^{1}(\mathbb{R})$ for some subsequence. Hence,
$$
\liminf_{n \to +\infty}\int_{\Re}(|u_{n}'|^{2}+|u_n|^{2}) = \int_{\Re}(|u'|^{2}+|u|^{2}),
$$
$$
\lim_{n \to +\infty} \int_{\Re}(u_n'u'+u_n u)= \int_{\Re}(|u'|^{2}+|u|^{2}),
$$
and
$$
\lim_{n \to +\infty}\int_{\Re}g(x,u_{n})(u-u_{n})=0.
$$
The above limits combined with (\ref{Z4}) yield
$$
\lambda_{n} \int_{\Re} V|u_{n}|^{2} \to 0.
$$
To prove that $(P_{{\cal O}})$ holds, we begin recalling that for all $v \in \Ka$,}
$$
\int_{\Re}[u_n'(v-u_n)'+(1+\lambda_{n} V(x))u_n(v-u_n)]\geq \int_{\Re}g(x,u_n)(v-u_n).
$$
Choosing $v \in  \widehat{\Ka}$, we get
$$
\int_{\Re}[u_n'(v-u_n)'+ u_n(v-u_n) - \lambda_n V(x)|u_n|^{2}]\geq \int_{\Re}g(x,u_n)(v-u_n).
$$
Taking the limit of $n \to \infty$ and using the Lemma \ref{convergencia} and (\ref{W1}), we derive that
$$
\int_{{\cal O}}[u'(v-u)' + u(v-u)] \geq \int_{{\cal O}}f(u)(v-u) \,\,\, \forall v \in \widehat{\Ka},
$$
finishing the proof. \fim

\begin{cor} \label{Linfinito} The sequences $(u_n)$ and $(w_n)$ satisfy the following limits
$$
\|w_n \|_{L^{\infty}({\overline{{\cal O}}}^{c})}, \|u_n \|_{L^{\infty}({\overline{{\cal O}}}^{c})} \to 0 \,\,\, \mbox{as} \,\,\, n \to +\infty.
$$

\end{cor}

\dem These limits are an immediate consequence of the continuous embedding $H^{1}({\overline{\Omega}}^{c}) \hookrightarrow
L^{\infty}({\overline{{\cal O}}}^{c})$ together with the limits $u_n \to u$ and $w_n \to w$ in $H^{1}(\Re)$ and of the fact that $u=w=0$ in ${\cal O}^{c}$. \fim

\vspace{0.5 cm}

\noindent \textbf{Proof of Theorem \ref{T1}}

The study made in this section allows us to prove that (\ref{Condicao}) holds for $\lambda$ large enough. We will show only (\ref{Condicao}) to $(u_n)$, because the argument is the same for $(w_n)$. Arguing by contradiction, we assume that there is $\lambda_n \to +\infty$ such that
\begin{equation}\label{W2}
\|u_n \|_{L^{\infty}({\Omega}^{c})} > a \,\,\,\,\, \forall n \in \mathbb{N}.
\end{equation}
From Lemma \ref{convergencia}, there is a subsequence of $(u_n)$, still denoted by itself, and $u \in H^{1}_{0}({\cal  O})$ such that
$$
u_n \to u \,\,\, \mbox{in} \,\,\, H^{1}(\Re).
$$
By Corollary \ref{Linfinito}, the below limit holds
$$
\|u_n \|_{L^{\infty}({\overline{{\cal O}}}^{c})} \to 0 \,\,\, \mbox{as} \,\,\, n \to +\infty,
$$
which implies that there is $n_{0} \in \mathbb{N}$ such that
$$
\|u_n \|_{L^{\infty}({\Omega}^{c})} < \frac{a}{2} \,\,\,\,\, \forall n \geq n_{0},
$$
obtaining a contradiction with  (\ref{W2}). This way, it follows that there is $\lambda^{*}>0$ such that the solution $u_\lambda$ satisfies
$$
u_\lambda (x) \leq a \,\,\, \forall x \in {\Omega}^{c} \,\,\, \mbox{and} \,\,\,  \lambda \geq \lambda^{*}.
$$
Now, by Remark \ref{R1}, we can conclude that $u_\lambda$ is a positive solution for $(P_\lambda)$ for $\lambda \geq \lambda^{*}$. \fim

\vspace{0.5 cm}

\noindent {\bf Acknowledgement}. Part of this work was done while the second author was visiting the Faculdade de Matem\'atica of the Universidade Federal do Par\'a (Brazil). In particular, he would like to express his deep gratitude to Prof. Giovany Figueiredo (UFPA) for his kind hospitality. The authors would like to thank to the referee for his/her suggestions to improve this paper.

\end{document}